\documentclass[12pt]{amsart}
\usepackage{amsmath,amscd,amssymb,amsfonts}
\newcommand{\h}{\hbox}
\newcommand{\q}{\quad}
\newcommand{\nin}{\noindent}
\newcommand{\bs}{\par\bigskip}
\newcommand{\ms}{\par\medskip}
\newcommand{\sk}{\par\smallskip}
\newcommand{\mtim}{\h{$\times$}}
\newcommand{\msum}{\h{$\sum$}}
\newcommand{\mopl}{\h{$\bigoplus$}}
\newcommand{\ssc}{\,\raise.15ex\h{${\scriptstyle\circ}$}\,}
\newcommand{\al}{\alpha}
\newcommand{\be}{\beta}
\newcommand{\de}{\delta}
\newcommand{\ko}{\overline{k}}
\newcommand\CH{{\rm CH}}
\newcommand{\Ga}{\Gamma}
\newcommand{\ga}{\gamma}
\newcommand{\Cor}{{\rm Cor}}

\newcommand{\Q}{{\mathbf Q}}
\newcommand{\R}{{\mathbf R}}
\newcommand{\RR}{{\mathcal R}}
\newcommand{\Z}{{\mathbf Z}}
\newcommand{\M}{{\mathcal M}}
\newcommand{\ges}{\geqslant}
\newcommand{\les}{\leqslant}
\begin{document}
\title{Ambiguity of certain Chow-K\"unneth projectors}
\author{Morihiko Saito}
\address{RIMS Kyoto University, Kyoto 606-8502 Japan}
\begin{abstract}
We show that the ambiguity of Murre's Chow-K\"unneth projector for degree 1
has certain good properties, assuming only that it factors through a Chow motive of a smooth irreducible curve. This is compatible with a picture obtained by using Beilinson's conjectural mixed motives, and implies for instance the independence of the Chow motive defined by the middle Chow-K\"unneth projector for surfaces.
\end{abstract}
\maketitle
\centerline{\bf Introduction}
\bs\nin
Let $X$ be a smooth projective variety of dimension $n$ over a field $k$ (assumed absolutely irreducible). Murre's conjecture (see [Mu1], [Mu2]) predicts the existence of Chow-K\"unneth projectors $\pi_0,\dots,\pi_{2n}$ in the ring of correspondences
$$\RR_X:=\Cor^0(X,X)=\CH^n(X\mtim X)_{\Q},$$
which are mutually orthogonal, and satisfy the following conditions:
\sk\nin
(1) $\sum_j\pi_j=[\Delta_X]$ in $\RR_X$, where $[\Delta_X]$ is the class of the diagonal.
\sk\nin
(2) The action of $\pi_i$ on $H^j(X_{\ko},\Q_l)$ is the identity for $j=i$, and vanishes otherwise. 
\sk
Here $H^j(X_{\ko},\Q_l)$ means the \'etale cohomology (with $l\ne{\rm char}\,k$) of the base change $X_{\ko}$ of $X$ by an algebraic closure $\ko$ of $k$.
\sk
For $\xi\in\RR_X$, we have the decompositions
$$\aligned&\q\xi=\sum_i{}_{(i)}\xi=\sum_j\xi_{(j)}=\sum_{i,j}{}_{(i)}\xi_{(j)}\q\h{with}\\
&{}_{(i)}\xi:=\pi_i\ssc\xi,\q\xi_{(j)}:=\xi\ssc\pi_j,\q{}_{(i)}\xi_{(j)}:=\pi_i\ssc\xi\ssc\pi_j,\endaligned
\leqno(3)$$
and these are respectively called the {\it left $\pi_i$-factor, right $\pi_j$-factor}, and $(\pi_i,\pi_j)$-{\it factor} of $\xi$.
\sk 
For $i=1$, the projector $\pi_1$ was constructed in loc.~cit.\ in such a way that $\pi_1$ factorizes through a Chow motive defined by an absolutely irreducible smooth projective curve $C$, i.e.,
$$\pi_1=\be\ssc\al\q\h{with}\q\al\in\Cor^0(X,C),\,\,\be\in\Cor^0(C,X),
\leqno(4)$$
where we assume $n\ges 2$.
This seems apparently weaker than the conditions coming from the construction of the projectors $\pi_1$ in loc.~cit. (Here it seems rather difficult to remove the absolute irreducibility condition on the curve $C$, see Remark~(1.6) below.) We may respectively replace $\al$ and $\be$ with $\al\ssc\pi_1$ and $\pi_1\ssc\be$ in (4). We may hence assume that the action of $\al$, $\be$ between $H^j(X_{\ko},\Q_l)$ and $H^j(C_{\ko},\Q_l)$ vanishes for $j\ne 1$.
\sk
Let $\zeta\in\CH_0(X)_{\Q}$ having degree 1 and corresponding to $\pi_0$, i.e.
$$\pi_0=\zeta\mtim[X]=\Ga_{[X]}\ssc{}^t\Ga_{\zeta},
\leqno(5)$$
where $\zeta$ is identified with $\Ga_{\zeta}\in\Cor^n(pt,X)$, and $\Ga_{[X]}\in\Cor^0(pt,X)$ is defined by $[X]$. Here $^t\Ga_{\zeta}$ denotes the transpose of $\Ga_{\zeta}$. (The factorization in (5) may be useful for non-experts in order to make some arguments easier to understand.) By hypothesis the projectors are mutually orthogonal, and we get
$$\pi_0\ssc\pi_1=0,\q\pi_1\ssc\pi_0=0.
\leqno(6)$$
Note that the last equality in (6) follows from (2) since $\pi_1\ssc\Ga_{[X]}\in\CH^0(X)_{\Q}$, see (1.1.2) below. It is unclear whether for any projector $\pi_1$ in $\RR_X$ satisfying only condition (2) for $i=1$, there is $\pi_0$ satisfying (5) and (6), see Remark~(1.6)(i) below.
\sk
Let $\pi'_1$ be another Chow-K\"unneth projector in $\RR_X$. Here we assume that $\pi'_1$ satisfies condition~(2) for $i=1$, and also (6) for $\pi'_0$ defined as in (5) with $\zeta$ replaced by $\zeta'$. Set
$$\de_i:=\pi'_i-\pi_i\q(i=0,1).$$
We have the following.
\ms\nin
{\bf Theorem~1.} {\it Assume $\pi_0=\pi'_0$, i.e.\ $\zeta=\zeta'$. Then $\de_1$ is orthogonal to $\pi_0$, and we have}
$$\de_1\ssc\pi_1=0,\q\q\pi'_1\ssc\de_1=\de_1.$$
\ms
The last two equalities respectively mean that the right $\pi_1$-factor $(\de_1)_{(1)}$ of $\de_1$ vanishes, and $\de_1$ coincides with its left $\pi_1$-factor $_{(1)}(\de_1)$ if we assume that $\pi_0=\pi'_0$ and condition~(4) is satisfied also for $\pi'_1$ (so that we can exchange $\pi_1$ and $\pi'_1$). Theorem~1 is reduced to the following.
\ms\nin
{\bf Proposition~1.} {\it Assume $\pi_0=\pi'_0$, i.e.\ $\zeta=\zeta'$. Then}
$$\pi'_1\ssc\pi_1=\pi_1.$$
\ms
Indeed, the first assertion of Theorem~1 on the orthogonality follows from (6), and the last equality of Proposition~1 is equivalent to either of the last two equalities of Theorem~1, since $\pi_1$ and $\pi'_1$ are idempotents.
\sk
In general, set
$$\de'_1:=\pi_0\ssc\pi'_1,\q\de''_1:=\de_1-\de'_1,\q\h{i.e.}\q\de_1=\de_1'+\de''_1.$$
Then
$$\de'_1=\Ga_{[X]}\ssc{}^t\Ga_{\eta}\q\h{with}\q\eta:=(^t\pi'_1)_*\,\zeta,\q\h{i.e.}\q{}^t\Ga_{\eta}={}^t\Ga_{\xi}\ssc\pi'_1.
\leqno(7)$$
Note that $\eta$ has degree 0 by condition (2) for $i=1$.
\ms\nin
{\bf Theorem~2.} {\it Assume condition~$(4)$ is satisfied also for $\pi'_1$. In the notation of $(3)$ we have}
$$\de'_1={}_{(0)}(\de'_1),\q{}_{(i)}(\de'_1)=0\,\,\,(i\ne 0),\q(\de'_1)_{(0)}=0,
\leqno(8)$$
$$\de''_1={}_{(1)}(\de''_1),\q{}_{(i)}(\de''_1)=0\,\,\,(i\ne 1),\q(\de''_1)_{(0)}=(\de''_1)_{(1)}=0.
\leqno(9)$$
\ms
Note that the first equalities of (8) and (9) respectively imply their second equalities (by using the orthogonality among the $\pi_i$). We do not have any assertion on $(\de'_1)_{(1)}$.
These assertions are compatible with a picture obtained by using Beilinson's conjectural mixed motives [Be] (see (1.7) below). By Theorem~2, we get the following.
\ms\nin
{\bf Corollary~1.} {\it Set $\pi_{[a,b]}:=\sum_{i=a}^b\pi_i$ for $a\les b$. The Chow motives defined by the projectors $\pi_0$, $\pi_1$, $\pi_{[2,2n-2]}$, $\pi_{2n-1}$, $\pi_{2n}$ do not depend on the choice of the Chow-K\"unneth projectors $\pi_0,\dots,\pi_{2n}$ such that condition~$(4)$ is satisfied for $\pi_1$ and $^t\pi_{2n-1}$. More precisely, for another set of Chow-K\"unneth projectors $\pi'_0,\dots,\pi'_{2n}$ satisfying the same conditions, we have
$$\pi_j\ssc\pi'_j\ssc\pi_j=\pi_j,\q\pi'_j\ssc\pi_j\ssc\pi'_j=\pi'_j,
\leqno(10)$$
for $j=0,1,2n-1,2n$ and also for $j=[2,2n-2]$, so that the Chow motives defined by $\pi_j$ and $\pi'_j$ are canonically isomorphic for each $j$ as above.}
\ms
For Chow motives, see (1.1) below.
Except for the case $j=[2,2n-2]$, this corollary is known to specialists, and can be proved by using the relation with the category of Abelian varieties in case $j=1$ or $2n-1$ as in [MNP], Th.~2.7.2.
\sk
We also have the following corollaries generalizing well-known results of Murre [Mu1], [Mu2] which were shown for $\pi_1,\pi_{2n-1}$ constructed in loc.~cit.
\ms\nin
{\bf Corollary~2.} {\it Assume $\pi_1$ and $^t\pi_{2n-1}$ satisfy condition~$(4)$. Then
$$(\pi_1)_*\CH^j(X)_{\Q}=0\q(j\ne 1),\q(\pi_{2n-1})_*\CH^j(X)_{\Q}=0\q(j\ne n),$$
and
$$\aligned{\rm Im}\bigl((\pi_1)_*:\CH^1(X)_{\Q}\to\CH^1(X)_{\Q}\bigr)&=\CH^1_{\rm hom}(X)_{\Q},\\
{\rm Ker}\bigl((\pi_{2n-1})_*:\CH^n(X)^0_{\Q}\to\CH^n(X)^0_{\Q}\bigr)&=T(X)_{\Q},\endaligned$$
where $\CH^1_{\rm hom}(X)_{\Q}$ is the subgroup of $\CH^1(X)_{\Q}$ consisting of cycles which are homologically equivalent to zero, $\CH^n(X)^0_{\Q}$ is the subgroup of $\CH^n(X)_{\Q}$ consisting of zero-cycles of degree $0$, and $T(X)_{\Q}\subset\CH^n(X)^0_{\Q}$ is the Albanese kernel.}
\ms\nin
{\bf Corollary~3.} {\it For $n=2$, the associated filtration on the Chow groups $($see $(1.3)$ below$)$ is independent of the Chow-K\"unneth projectors $\pi_i$ satisfying the hypothesis in Corollary~$2$.}
\ms
These can also be proved directly by generalizing Murre's argument [Mu1], [Mu2], and using a commutative diagram as in [MNP], p.~75 (see Remark~(2.6)(i) below).
For the relation with the conjectural filtration of Bloch [Bl] and Beilinson [Be] on the Chow groups, see [Ja]. It seems rather difficult to remove the hypothesis in the corollaries, see Remark~(1.4) below. It does not seem trivial to generalize Corollaries~1 and 2 to the situation treated in [Sa2]. 
\sk
This work is partially supported by Kakenhi 24540039.
I thank J.~Murre and J.~Nagel for useful discussions about Chow-K\"unneth decompositions.
\ms
In Section 1 we review some facts from Chow motives.
In Section 2 we prove the main theorems and their corollaries.
\bs\bs
\centerline{\bf 1. Chow motives}
\bs\nin
In this section we review some facts from Chow motives.
\ms\nin
{\bf 1.1.~Correspondences.} Let $X,Y$ be smooth proper algebraic varieties $X,Y$ defined over a field $k$. Assume $X$ purely $n$-dimensional. The group of correspondences with rational coefficients is defined by
$$\Cor^i(X,Y)=\CH^{n+i}(X\mtim Y)_{\Q}. $$
The composition of correspondences is defined as in the case of Grothendieck motives (see [Kl], [Ma]), and is denoted by
$$\eta\ssc\xi\in\Cor^{i+j}(X,Z)\,\,\,\h{for}\,\,\,\xi\in\Cor^i(X,Y),\,\,\eta\in\Cor^{j}(Y,Z).
\leqno(1.1.1)$$
If $X=pt$ so that
$$\Cor^i(pt,Y)=\CH^i(Y)_{\Q},$$
then the composition is identified with the action of $\eta$ on the Chow groups
$$\eta_*:\CH^j(Y)_{\Q}\to\CH^{i+j}(Z)_{\Q}.$$
If $i=j=0$, $X=pt$, and $Y=Z$ is absolutely irreducible in (1.1.1), then we have
$$\eta\ssc\xi=c_0\,\xi,
\leqno(1.1.2)$$
where $c_0$ is the (unique) eigenvalue of the action of $\eta_*$ on $H^0(Y_{\ko},\Q_l)$.
\sk
For $\xi\in\CH^i(Y)_{\Q}$, the correspondence defined by it is sometimes denoted by
$$\Ga_{\xi}\in\Cor^i(pt,Y)_{\Q}\,\bigl(=\CH^i(Y)_{\Q}\bigr).$$
If $i=\dim Y$, we have $^t\Ga_{\xi}\in\Cor^0(Y,pt)$, and
$$^t\Ga_{\xi}\ssc\Ga_{[Y]}=\deg\xi\in\Cor^0(pt,pt)=\Q.
\leqno(1.1.3)$$
\sk
For a projector $\ga\in\Cor^0(X,X)$ and $m\in\Z$, we have the associated Chow motive $(X,\ga,m)$ so that for another Chow motive $(X',\ga',m')$, we have by definition (see [MNP])
$${\rm Hom}\bigl((X,\ga,m),(X',\ga',m')\bigr)=\ga'\ssc\Cor^{m'-m}(X,X')\ssc\ga.
\leqno(1.1.4)$$
In case $m=0$, $(X,\ga,0)$ is denoted by $(X,\ga)$, and is called the Chow motive defined by the projector $\ga$. The Tate twist $(p)$ for $p\in\Z$ is defined so that
$$(X,\ga,m)(p)=(X,\ga,m+p).
\leqno(1.1.5)$$

\ms\nin
{\bf 1.2.~Codimension 1 case.}
Assume $X,Y$ absolutely irreducible. It is well-known (see e.g.\ [Mu1], [MNP], [Sch]) that there is a canonical isomorphism due to Grothendieck
$$\frac{\CH^1(X\mtim Y)_{\Q}}{{pr}_1^*\CH^1(X)_{\Q}+{pr}_{2}^*\CH^1(Y)_{\Q}}={\rm Hom}({\rm Alb}_X,{\rm Pic}^0_Y)_{\Q},
\leqno(1.2.1)$$
together with an injection
$${\rm Hom}\bigl({\rm Alb}_X,{\rm Pic}^0_Y)_{\Q}\,\hookrightarrow\,{\rm Hom}(H^{2n-1}(X_{\ko},\Q_l)(n-1),H^1(Y_{\ko},\Q_l)\bigr).
\leqno(1.2.2) $$
\sk
Let $\xi\in\CH_{0}(X)_{\Q}$, $\eta\in\CH_{0}(Y)_{\Q}$ with degree one.
Then the left-hand side of (1.2.1) is isomorphic to
$$\bigl\{\ga\in\Cor^{1-n}(X, Y)_{\Q}\,\,\big|\,\,{}^t\Ga_{\eta}\ssc\ga=0,\,\,\ga\ssc\Ga_{\xi}=0\bigr\}.
\leqno(1.2.3)$$
We say that $\ga$ is {\it left-normalized by} $\eta$ if $^t\Ga_{\eta}\ssc\ga=0$, i.e.\ if $(^t\ga)_*\eta=0$, and {\it right-normalized by} $\xi$ if $\ga\ssc\Ga_{\xi}=0$, i.e.\ if $\ga_*\xi=0$.
\ms\nin
{\bf 1.3.~Filtration on Chow groups.} It is conjectured in [Mu2] that the action of the Chow-K\"unneth projector $\pi_i$ on $\CH^j(X)_{\Q}$ vanishes unless $i\in[j,2j]$. The filtration $F$ on $\CH^j(X)_{\Q}$ is defined in loc.~cit. by
$$F^k\CH^j(X)_{\Q}:=\bigcap_{i>2j-k}{\rm Ker}\,(\pi_i)_*=\sum_{i\le 2j-k}{\rm Im}\,(\pi_i)_*\subset\CH^j(X)_{\Q}.
\leqno(1.3.1)$$
\ms\nin
{\bf 1.4.~Remark.} The hypothesis in Corollary~1 seems rather essential, since it is related with a problem of phantom motives, i.e. nontrivial motives such that the action of the defining idempotent on the \'etale cohomology vanishes. (Indeed, in case such a motive exists, we cannot exclude the possibility that $\pi_1$ contains it.) It is conjectured that such a motive does not exist. This problem is rather nontrivial since it is closely related with Bloch's conjecture [BKL], see [GP], [KMP], [Sa1]. Here we can construct mutually orthogonal projectors $\pi_0,\dots,\pi_4$ quite explicitly for a surface with $p_g=0$, but it is unclear whether condition~(1) is satisfied, and the difference of both sides in condition~(1) might define a phantom motive.
\ms
The following lemma will be used in the proof of Proposition~1 in (2.1) below.
\ms\nin
{\bf 1.5.~Lemma.} {\it Let $\be$ be as in $(4)$. Assume the action
$$\be_*:H^j(C_{\ko},\Q_l)\to H^j(X_{\ko},\Q_l)
\leqno(1.5.1)$$
vanishes for $j=2$. Then there is $\xi\in\CH_0(C)_{\Q}$ of degree $1$ whose image in $\CH^1(X)_{\Q}$ by $\be_*$ vanishes, i.e. $\beta$ is right-normalized by $\xi$ in the sense of $(1.2)$.}
\ms\nin
{\it Proof.}
The action (1.5.1) is surjective for $j=1$ by using the factorization (4) together with condition~(2). It implies the surjectivity of the morphism
$$\be_*:\CH_{\rm hom}^1(C)_{\Q}\to\CH_{\rm hom}^1(X)_{\Q},$$
where $\CH_{\rm hom}^1(X)_{\Q}$ denotes the subgroup consisting of homologically equivalent to zero cycles. So the assertion follows.
\ms\nin
{\bf 1.6.~Remarks.}
(i) Let $\pi_1$ be a projector satisfying condition (2) for $i=1$. It is unclear whether there is $\pi_0$ satisfying (5) and the first equality of (6) for this $\pi_1$.
In fact, the argument in the proof of Lemma~(1.5) cannot be generalized to this case, since the action of $^t\pi_1$ on the Albanese kernel may be non-surjective.
\ms
(ii) Lemma~(1.5) implies that $\pi_1$ factors through a motive defined by a Chow-K\"unneth projector $\pi^C_1$ for the curve $C$, where the projector is defined by
$$\pi^C_1:=[\Delta_C]-\Ga_{[C]}\ssc{}^t\Ga_{\eta}-\Ga_{\eta'}\ssc{}^t\Ga_{[C]},$$
with $\eta,\eta'\in\CH_0(C)_{\Q}$ of degree 1. Here $\eta'$ is the one given by Lemma~(1.5), but $\eta$ may be arbitrary. This factorization is equivalent to
$$\be\ssc\Ga_{[C]}\ssc{}^t\Ga_{\eta}\ssc\al+
\be\ssc\Ga_{\eta'}\ssc{}^t\Ga_{[C]}\ssc\al=0,$$
and is reduced to
$$\be\ssc\Ga_{[C]}=0,\q\be\ssc\Ga_{\eta'}=0.$$
Here the first vanishing is easy, and the second follows from Lemma~(1.5).

\ms\nin
{\bf 1.7.~Beilinson's conjectural mixed motives.}
It has been conjectured by Beilinson [Be] that there is an abelian category of mixed motives $\M$ for varieties over $k$ together with $\R\Ga_{\M}(X)\in D^b\M$ for smooth projective varieties $X$ over $k$, such that we have a decomposition
$$\R\Ga_{\M}(X)\cong\mopl_{j=0}^{2n}\,H^j_{\M}(X)[-j],
\leqno(1.7.1)$$
(with $n:=\dim X$) and also a canonical isomorphism
$$\CH^j(X)_{\Q}={\rm Hom}_{D^b\M}\bigl(\Q_{\M},\R\Ga_{\M}(X)(j)[2j]\bigr),
\leqno(1.7.2)$$
where $\Q_{\M}=\R\Ga_{\M}({\rm Spec}\,k)\in\M$, and $(j)$ is the Tate twist. Then the canonical truncations $\tau_{\les k}$ (see [De]) on the complex $\R\Ga_{\M}(X(j)[2j]$ define a filtration on $\CH^j(X)_{\Q}$, which is called the (conjectural) filtration of Bloch and Beilinson.
It is further conjectured that we have an isomorphism
$$\RR_X={\rm Hom}_{D^b\M}\bigl(\R\Ga_{\M}(X),\R\Ga_{\M}(X)\bigr),
\leqno(1.7.3)$$
so that a decomposition (1.7.1) should correspond to a set of Chow-K\"unneth projectors $\pi_0,\dots,\pi_{2n}$ (see also [Ja]).
By (1.7.1) and (1.7.3), any $\xi\in\RR_X$ is expressed by
$$\xi_{i,j}\in{\rm Ext}^{j-i}_{\M}\bigl(H_{\M}^j(X),H_{\M}^i(X)\bigr)\q\h{for}\,\,\,i\les j,$$
If $\xi$ is an idempotent, then these should satisfy the relations
$$\xi_{i,j}=\msum_{i\les k\les j}\,\xi_{i,k}\ssc\xi_{k,j}\q\h{for}\,\,\,i\les j.$$
Assume furthermore $\xi_{1,1}=id$, and $\xi_{i,i}=0$ for $i\ne 1$. Then we get by induction on $j-i$
$$\xi_{i,j}=0\q\h{if}\,\,\,i>1,
\leqno(1.7.4)$$
and moreover
$$\xi_{0,j}=\xi_{0,1}\ssc\xi_{1,j}\q\h{for}\,\,\,j>1.
\leqno(1.7.5)$$
Note that $\xi_{i,j}$ should correspond in the notation of (3) to
$${}_{(i)}\xi_{(j)}.$$
We can apply these to the case $\xi=\pi'_1$.
\ms
The following will be used in the proof of Corollary~2.
\ms\nin
{\bf 1.8.~Lemma.} {\it Let $\varphi,\delta$ be two endomorphisms of a vector space $V$. Then}
$$\aligned{\rm Im}\,(\varphi+\delta)\subset {\rm Im}\,\varphi&\q\h{if}\q\de=\varphi\ssc\de,\\{\rm Ker}\,\varphi\subset {\rm Ker}(\varphi+\delta)&\q\h{if}\q\de=\de\ssc\varphi.\endaligned$$

\bs\bs
\centerline{\bf 2. Proof of the main theorems}
\bs\nin
In this section we prove the main theorems and their corollaries.
\ms\nin
{\bf 2.1.~Proof of Proposition~1.} Set
$$\be'':=\pi'_1\ssc\be.$$
We have to show
$$\be=\be''\q\h{in}\q\Cor^0(C,X)=\CH^1(C\mtim X)_{\Q},
\leqno(2.1.1)$$
by replacing $\al,\be$ appropriately (without changing $\pi_1$).
\sk
By (1.1.3) together with (5), the first condition of (6) is equivalent to
$$^t\Ga_{\xi}\ssc\pi_1=0.
\leqno(2.1.2)$$
Using the hypothesis $\xi=\xi'$, we then get
$$^t\Ga_{\xi}\ssc\be=0,\q^t\Ga_{\xi}\ssc\be''=0,
\leqno(2.1.3)$$
since we may replace $\be$ with $\pi_1\ssc\be$.
\sk
On the other hand, we have for some $\eta\in\CH_0(C)_{\Q}$ with degree 1
$$\be\ssc\Ga_{\eta}=0,\q\h{hence}\q\be''\ssc\Ga_{\eta}=0.
\leqno(2.1.4)$$
Indeed, this is equivalent to the assertion that there is a zero-cycle of degree 1 belonging to the kernel of
$$\be_*:\CH^1(C)_{\Q}\to\CH^1(X)_{\Q},$$
and follows from Lemma~(1.5).
\sk
By (1.2.1--3) and (2.1.3--4), the assertion (2.1.1) is thus reduced to
$$\be_*=\be''_*:H^1(C_{\ko},\Q_l)\to H^1(X_{\ko},\Q_l),
\leqno(2.1.5)$$
and follows from the assumption that the action of $\pi'_1$ on $H^1(X_{\ko},\Q_l)$ is the identity. This finishes the proof of Proposition~1.
\ms\nin
{\bf 2.2.~Remarks.} (i) The dual assertion holds for $\pi_{2n-1}$.
\ms
(ii) We have a similar assertion for $\pi_0$ and $\pi_{2n}$. Indeed, $\pi_0$ is defined by (5) depending on $\xi\in\CH_0(X)_{\Q}$ with degree 1, and, setting $\de_0:=\pi'_0-\pi_0$, we have
$$\pi'_0\ssc\pi_0=\pi_0,\q\pi'_0\ssc\de_0=\de_0,
\leqno(2.2.1)$$
if $\pi'_0$ is associated to $\xi'\in\CH_0(X)_{\Q}$ with degree 1.
\ms\nin
{\bf 2.3.~Proof of Theorem~2.} The equalities in (8) follow from (7) together with a remark after (6). For the proof of (9), set
$$\pi''_1:=\pi'_1-\de'_1=(1-\pi_0)\ssc\pi'_1.$$
Since $\pi_0$ is an idempotent, we have
$$\pi_0\ssc\pi''_1=0.$$
Moreover, we have
$$\pi''_1\ssc\pi_0=\pi'_1\ssc\pi_0=0,$$
by the same argument as in a remark after (6). So $\pi''_1$ is an idempotent, since
$$\pi''_1{}^2=(1-\pi_0)\ssc\pi'_1\ssc(1-\pi_0)\ssc\pi'_1.$$
Then (9) follows from Theorem~1 applied to $\pi''_1$ and $\pi_1$, since we have by definition
$$\de''_1=\pi''_1-\pi_1.$$
This finishes the proof of Theorem~2.
\ms\nin
{\bf 2.4.~Proof of Corollary~1.} For $j=0,1,2n-1,2n$ and also for $j=[2,2n-2]$, we have to show
$$\pi_j\ssc\pi'_j\ssc\pi_j=\pi_j.$$
Set
$$\de_j:=\pi'_j-\pi_j.$$
Since $\pi_j\ssc\pi_j\ssc\pi_j=\pi_j$, the assertion is equivalent to
$$\pi_j\ssc\de_j\ssc\pi_j=0.
\leqno(2.4.1)$$
For $j=1$, this follows from the second equality of (8) and the last one of (9) in Theorem~2. For $j=2n-1$, we can apply these to $^t\pi_{2n-1}$.
The argument is similar for $j=0$ and $2n$ by using (2.2.1).
\sk
So it remains to show the case $j=[2,2m-2]$, where we have by definition
$$\pi'_j-\pi_j=-\sum_{j\in\{0,\,1,\,2n-1,\,2n\}}\de_i.$$
It is then enough to prove for $i=0,1,2n-1,2n$
$$\biggl(1-\sum_{a\in\{0,\,1,\,2n-1,\,2n\}}\pi_a\biggr)\ssc\de_i\ssc\biggl(1-\sum_{a\in\{0,\,1,\,2n-1,\,2n\}}\pi_a\biggr)=0.
\leqno(2.4.2)$$
For $i=1$, it is sufficient to show (2.4.2) with $\de_1$ replaced by $\de'_1$ and also by $\de''_1$ in the notation of Theorem~2. Then the assertion follows from the first and second equalities of (8) and (9). For $i=2n-1$, we can apply the same argument to $^t\pi_{2n-1}$.
The argument is similar for $i=0$ or $2n$ by using (2.2.1).
This finishes the proof of Corollary~1.

\ms\nin
{\bf 2.5.~Proof of Corollary~2.} The first assertions easily follows from the factorization~(4) and its transpose by considering the factorization of the action on the Chow groups. As for the remaining assertions, the independence of the left-hand side follows from Lemma~(1.8) by using the first equalities of (8) and (9) in Theorem~2 applied to $\pi_1$, $\pi'_1$, and also to $^t\pi_{2n-1}$, $^t\pi'_{2n-1}$. (The details are left to the reader.) So the assertions are reduced to the result of Murre, see [Mu1], [Mu2], [MNP].
This finishes the proof of Corollary~2. (See also Remark~(2.6)(i) below.)
\ms\nin
{\bf 2.6.~Remarks.}
(i) We can also prove Corollary~2 directly by generalizing Murre's argument in [Mu1], [Mu2] using the action of $\pi_1$ and $\pi_{2n-1}$ on the Picard and Albanese varieties together with a commutative diagram as in [MNP], p.~75. Here the irreducibility of $C$ seems to be unnecessary. (The details are left to the reader.)
\ms
(ii) It is unclear whether the main theorems in this paper can be generalized to the case where the curve $C$ in condition~(4) is not necessarily irreducible, since the same argument does not seem to work. It is also unclear whether the argument in Remark~(1.6)(ii) can be generalized to this case.

\end{document}